\documentclass[12pt,twoside]{amsart}
\usepackage{amsmath}
\usepackage{amsfonts}
\usepackage{amssymb}
\usepackage{amsthm}
\usepackage[english]{babel}

\newtheorem{theo}{Theorem}[section]
\newtheorem{lem}[theo]{Lemma}
\newtheorem{prop}[theo]{Proposition}
\newtheorem{cor}[theo]{Corollary}

\newcommand{\E}{\ensuremath{\mathbb {E}}}
\renewcommand{\P}{\ensuremath{\mathbb {P}}}

\newcommand{\N}{{\mathcal N}}

\newcommand{\F}{{\mathcal F}}

\newcommand{\B}{{\cal B}}

\renewcommand{\H}{{\mathcal H}}

\def\E{\mathbb{E}}

\def\B{{\cal B}}

\begin{document}


\title{ Pointwise ergodic theorems with rate and application to limit
theorems for stationary processes}

\medskip

 \author{ Christophe Cuny}
\curraddr{Equipe ERIM, University of New Caledonia, B.P. R4, 98800 Noum\'ea,
New Caledonia}
\email{cuny@univ-nc.nc}
\thanks{}



\medskip

\begin{abstract}
We obtain pointwise ergodic theorems with rate under conditions
expressed in terms of the convergence of  series involving $\|\sum_{k=1}
^nf\circ \theta^k\|_2$, improving previous results. Then, using known results
on martingale approximation, we obtain some LIL for stationary ergodic processes and
quenched central limit theorems for functional of Markov chains.
The proofs are based on the use of the spectral theorem and, on a recent
work of Zhao-Woodroofe extending a  method of Derriennic-Lin.
\end{abstract}

\maketitle

\section{Introduction}
One of the goals of the present paper is to obtain limit theorems
for $\{S_n:=\sum_{k=1}^nX_k\}$, where $\{X_n\}$ is a strictly stationary
process in $L^2$. In particular, we will be interested in the
Law of the Iterated Logarithm or,in the case where $\{X_n\}$ is given by
the functional of a Markov chain, in the quenched Central Limit Theorem.
\\
We will follow a classical line, using known approximations of
$\{S_n\}$ by a martingale $\{M_n\}$, yielding to a control of
$\{\|S_n-M_n\|_2\}$. \\
Then, it remains to obtain a.s. convergence results  on $\{S_n-M_n\}$
under the control of $\{\|S_n-M_n\|_2\}$. \\
This approach was used recently in few papers, see e.g. \cite{DL1},
\cite{DL2}, \cite{ZW}, \cite{CL} and \cite{Wu}.  \\
The first four papers used a martingale approximation
initiated by Kipnis-Varadhan \cite{KV} (and developped by Maxwell-Woodroofe
\cite{MW}) and results or extensions of Derriennic-Lin
\cite{DL} to obtain pointwise ergodic theorems for $\{S_n-M_n\}$,  while
Wu \cite{Wu} used another way to control $\{\|S_n-M_n\|_2\}$ and a
different way to obtain pointwise ergodic theorems (which is not efficient
to obtain LIL).
\\
In this paper we will show how the use of spectral tools may allow us
to obtain better pointwise ergodic theorems, using the approach of \cite{ZW}.
Then we use the martingale approximations in \cite{MW} or in \cite{Wu}
to obtain LIL and quenched central limit theorems, improving the results
of the previously mentionned papers. \\

Our paper is organised as follows.
In section 2, we collect some more or less known
 results of Gaposhkin about unitary operators (see \cite{Gap1}, \cite{Gap2}
and \cite{Gap3}) and extend them to normal operators. In section 3, we obtain
ergodic theorems with rates for  measure-preserving transformation and
functions in $L^2$. We focus our study to rates close to the critical rate
$\sqrt n$, improving former results of \cite{Wu}, \cite{ZW} and
\cite{CL} but also of \cite{Gap2}, \cite{CoL} and \cite{Weber}. In section 4,
we establish our LIL and quenched CLT. In section 5, we look at the particular
case of Markov chains with normal transition operator, which yields improved
results and a nice control of $\{\|S_n-M_n\|\}$. Finally, in section 6
we give examples allowing to compare our different results. The proofs of
technical nature are left to the appendix.

\section{Spectral criteria for the norm convergence of some series}

In this section, we give more or less known results about the norm convergence
of some power series associated with a normal operator.\\

In all the section $T$ will be a normal contraction of a Hilbert space $\H$
(i.e. $T^*T=TT^*)$ and $f\in \H$.
Denote $\mu_f$ the spectral measure of $f$, that is, $\mu_f$ is a finite
positive measure on the Borel sets of the closed unit disk $\overline D$,
such that for every $(a_0,\ldots ,a_n)\in {\bf C}^{n+1}$,
$\|\sum_{k=1}^n a_k T^kf\|^2=\int_{\overline D} |\sum_{k=1}^n a_k z^kf|^2
\mu_f(dz)$.\\

For every $n\ge 1$, write
$$D_n:=\{z=r{\rm e}^{2i\pi \theta}~:~1-\frac{1}{n}\le r\le 1,~-\frac{1}{n}\le
\theta \le \frac{1}{n}\}.$$
Hence $D_1=\overline D$ is the closed unit disk.\\
On can see (see e.g. \cite{CL}) that there exists $C>1$ such
 \begin{equation}\label{Dn}
\frac{n}{C}\le \frac{1}{|1-z|}\le C n,\qquad \forall z\in D_n-D_{n+1}.
\end{equation}

Write $U_n(f):=\sum_{k=1}^nT^k f$.
\begin{lem}\label{lem}
Let $\varphi ~:~D_1\rightarrow {\bf C}$ continuous on $D_1-\{1\}$ such that
there exists a non decreasing $\chi~:~{\bf R}^{+*}\rightarrow {\bf R}^{+*}$,
such that $|\varphi(z)|\underset{1}\sim \chi(\frac{1}{|1-z|})$. Moreover
assume that $\chi(x)/x^\alpha$ is non increasing for some $0<\alpha<2$ and
that there exists $\tau>1$ such that $\chi(2^{n+1})\ge \tau\chi(2^n)$.
Then, the  following are equivalent
\begin{itemize}
\item [(i)]   $\int_{D_1}|\varphi(z)|\mu_f(dz)<+\infty$.
\item [(ii)]  $\sum_{n\ge 1} \chi(2^n)\mu_f(D_{2^n})<+\infty$.
\item [(iii)] $\sum_{n\ge 1}\frac{\chi(n)}{n}\mu_f(D_n)<+\infty$.
\item [(iv)] $\sum_{n\ge 1} \frac{\chi(n)\|U_n(f)\|^2}{n^3}<+\infty$.
\end{itemize}
\end{lem}
{\bf Remarks :}  This result is a generalization of  \cite[Lemma 5]{Gap1}
stated for unitary operators. There is a continuous version of that lemma in
\cite[Lemma 1]{Gap2}.\\
{\bf Proof :}\\
We first prove $(i)\Leftrightarrow (ii)$.\\
We have $\int_{D_1}|\varphi(z)|\mu_f(dz)=\sum_{n\ge 1}\int_{D_{2^n}-
D_{2^{n+1}}}|\varphi(z)|\mu_f(dz)$. Using \eqref{Dn} and the monotony
assumptions on $\chi$, there exists $C>1$, such that
$$
\frac{1}{C}\sum_{n\ge 1}\chi(2^n)(\mu_f(D_{2^n})-\mu_f(D_{2^{n+1}})
\le \int_{D_1}|\varphi(z)|\mu_f(dz)\le C
\sum_{n\ge 1}\chi(2^n)(\mu_f(D_{2^n})-\mu_f(D_{2^{n+1}}),
$$
which proves $(ii)\Rightarrow (i)$. \\
On the other hand, using Fubini, since all the terms are non negative,
\begin{gather*}
\sum_{n\ge 1}(\chi(2^n)-\chi(1)(\mu_f(D_n)-\mu_f(D_{n+1})) =
\sum_{n\ge 1}(\mu_f(D_n)-\mu_f(D_{n+1})\sum_{k=0}^{-1}n(\chi(2^{k+1})-
\chi(2^k))\\
 =\sum_{k\ge 1}(\chi(2^{k+1})-\chi(2^k))\sum_{n\ge k}(\mu_f(D_n)-\mu_f(D_{n+1})
  =\sum_{k\ge 1}(\chi(2^{k+1})-\chi(2^k))\mu_f(D_k)\\
\ge (c-1)\sum_{k\ge 1} \chi(2^k)\mu_f(D_k)
\end{gather*}
which yields that $(i)\Rightarrow (ii)$.\\
For $(ii)\Leftrightarrow (iii)$ ), using $\{\chi(n)\}$ and $\{\mu_f(D_n)\}$
are monoton, we have
$$
\frac{1}{2}\sum_{k\ge 0}\chi(2^k)\mu_f(D_{2^{k+1}})\le
\sum_{k\ge 0} \sum_{n=2^k}^{2^{k+1}-1}\frac{\chi(n)\mu_f(D_n)}{n}
\le \sum_{k\ge 0}\chi(2^{k+1})\mu_f(D_{2^{k}}).
$$
which implies the desired equivalence since $\chi(x)/x^\alpha$ is non
increasing.\\
Let prove $(iii)\Leftrightarrow (iv)$.\\
By \cite{CL}, formulae $(10)$ and $(11)$, there exists $K>0$ such that for
every $n\ge 1$,
\begin{equation}\label{Sn}
\frac{n^2}{K}\mu_f(D_n)\le \|U_n(f)\|^2
\le \sum_{j=1}^{n-1}(2j+1)\mu_f(D_j).
\end{equation}
Hence $(iv)$ implies $(iii)$.\\
Assume that $\sum_{n\ge 1} \frac{\chi(n)\mu_f(D_n)}{n}<+\infty$.
Hence, using that  $\chi(x)/x^\alpha $ is non increasing, say for
$x\ge n_0\in{\bf N}$, we obtain, by \eqref{Sn},
\begin{gather*}
\sum_{n\ge n_0}\chi(n)\frac{\|U_n(f)\|^2}{n^3}\le
\sum_{n\ge n_0}\frac{\chi(n)}{n^3}
\sum_{j=1}^{n-1}3j\mu_f(D_j)\\
\le 3\sum_{j= 1}^{n_0-1}j\mu_f(D_j)\sum_{n\ge n_0}\frac{\chi(n)}
{n^\alpha}\frac{1}{n^{3-\alpha}}+
 3\sum_{j\ge n_0}j\mu_f(D_j)\sum_{n\ge j+1}\frac{\chi(n)}
{n^\alpha}\frac{1}{n^{3-\alpha}}\\
\le K+ K\sum_{j\ge n_0}\frac{\chi(j)\mu_f(D_j)}{j}.
\end{gather*}
\hfill $\square$

Following \cite{DL}, for every contraction $T$, we define the operator
$\sqrt{I-T}$ by $ \sqrt{I-T}:=\sum_{n\ge 0}\delta_n T^n$, where
$\sqrt{1-x}=\sum_{n\ge 0}\delta_n x^n$, where $\delta_n<0$ $\forall n\ge 1$
and $\sum_{n\ge 1}\delta_n<+\infty$.

\begin{prop}\label{sqrt}
The following are equivalent
\begin{itemize}
\item [(i)]  $f\in \sqrt {I-T}\H$.
\item [(ii)]  $\int_{D_1}\frac{\mu_f(dz)}{|1-z|}<+\infty$.
\item [(iii)] $\sum_{n\ge 1} \frac{\|U_n(f)\|^2}{n^2}<+\infty$.
\end{itemize}
\end{prop}
\noindent{\bf Proof :}\\
$(i)$ and $(ii)$ are equivalent, by  Theorem 4.4 of \cite{DL}. It can also be
deduced from Proposition \ref{CoL} (below) with $\psi(x):=
\sqrt x$, using that $(i)$ is equivalent to the norm convergence of
$\sum_{n\ge 0}c_nT^n$, where $\sum_{n\ge 0}c_nx^n=(1-x)^{-1/2}$, $0\le x<1$.
The equivalence of $(ii)$ and $(iii)$ follows from Lemma \ref{lem}, with
$\chi(x)=x$.
\hfill $\square$

\medskip

Define $ D:=  D(0,1)$.
\begin{prop}\label{CoL}
Let $\varphi (z)=\sum_{n\ge 0}a_n z^n$ be a power
series converging and continuous on $\overline D-\{1\}$, with $\{a_n\}$
 non negative.
 Assume that there exists a
non negative continuous function
$\psi~:~{\bf R}^{+*}\rightarrow {\bf R}^{+*}$, such that $|\varphi(z)|
\underset{1}\sim \psi(\frac{1}{|1-z|})$ ($\psi(+\infty)=+\infty$) and
$\sup_{m\ge 0} |\sum_{n=0}^ma_nz^n|\le C
\psi(\frac{1}{|1-z|})$ for every $z\in \overline D$. Then, for every normal
 contraction $T$ on $\H$ and any $f\in \H$ the following are equivalent
\begin{itemize}
\item [(i)]   $\sum_{n\ge 0} a_n T^n(f)$ converges in norm.
\item [(ii)]  $\sum_{n\ge 0} a_n T^n(f)$ converges weakly.
\item [(iii)] $\sup_{m\ge 0}\|\sum_{n=0}^ma_nT^n(f)\|<+\infty$.
\item [(iv)]  $\int_{\overline D} \psi^2 (\frac{1}{|1-z|})\mu_f(dz)<+\infty$.
\end{itemize}
\end{prop}
\noindent {\bf Remark:}\\
The equivalence of $(i)-(iv)$ was obtained in \cite{CoL} in the case treated
in Proposition \ref{Cohen} below.\\
{\bf Proof :}\\
$(i)\Rightarrow (ii)$ is clear and $(ii)\Rightarrow (iii)$ follows from
Banach-Steinhaus.\\
Let show $(iii)\Rightarrow (iv)$. By the spectral theorem,
for every $m\ge 0$, we have
$$
\|\sum_{n=0}^ma_nT^n(f)\|^2=\int_{\overline D}|\sum_{n=0}^m a_n
z^n|^2\mu_f(dz).
$$
Since $\{a_n\}$ are non negative and $\psi(+\infty)=+\infty$,
hence $\sum_{n\ge 0}a_n=+\infty$, and,
by $(iii)$, $\mu_f(\{1\})=0$. Hence, by assumption $\{\sum_{n=0}^m a_n
z^n\}$ converges for $\mu_f$-a.e. $z\in D_1$ to  $\varphi(z)$. By Fatou's
lemma we obtain
$$
\int_{\overline D}|\varphi(z)|^2\mu_f(dz)=\int_{\overline D}
\liminf_{m\rightarrow +\infty}|\sum_{n=0}^ma_nz^n|^2\mu_f(dz)\le
\sup_{m\ge 0}|\|\sum_{n=0}^ma_nT^n(f)\|^2<+\infty
$$
Thus $(iv)$ follows from the continuity of $\psi$ and $\varphi$, and
from the assumption
$|\varphi(z)|\underset{1}\sim \psi(\frac{1}{|1-z|})$.\\
Proof of $(iv)\Rightarrow (i)$. \\
Since $\psi(+\infty)=+\infty$, $(iv)$ implies that $\mu_f(\{1\})=0$.
Hence $\{\sum_{n=0}^m a_nz^n\}$ converges for $\mu_f$-a.e. $z$ to  $\varphi(z)$
and is dominated by $\psi(\frac{1}{|1-z|})$, so the dominated convergence
theorem implies that $\{\sum_{n=0}^m a_nz^n\}$ converges in $L^2(\mu_f)$,
which implies by the spectral theorem that $\{\sum_{n=0}^ma_nT^n(f)\}$ is
a Cauchy sequence, hence $(i)$. \hfill $\square$

There are plenty of power series for which the assumptions of Proposition
\ref{CoL} are satisfied. We give an application below. Our main applications
 will concern a family of particular power series introduced in \cite{ZW}.\\

\begin{prop}\label{Cohen}
The following are equivalent
\begin{itemize}
\item [(i)]   The series $\sum_{n\ge 1}\frac{T^nf}{n}$ converges in norm.
\item [(ii)]  $\int_{D_1}\log^2|1-z| \mu_f(dz) <+\infty$.
\item [(iii)] $\sum_{n\ge 1} \frac{\|U_n(f)\|^2\log n}{n^3}<+\infty$.
\end{itemize}
\end{prop}
\noindent {\bf Proof :}\\
$(i)$ and $(ii)$ are equivalent by \cite{CoL} (see also \cite{AL} for the case
of unitary or symmetric operators). One could apply the previous proposition
with $\varphi(z):=\sum_{n\ge 1}\frac{z^n}{n}$ and $\psi:=\log$. It is not hard
to see that the assumptions on $\varphi$ and $\psi$ are satified (it follows
from (\cite{Zygmund} Ch. I, p.2) and classical computations, it is also done
in \cite{CoL}). \\
To prove the equivalence to $(iii)$ one has to redo the proof of Lemma
\ref{lem} with $\psi:=\log^2$.
\hfill $\square$\\
{\bf Remarks:}\\
{\bf 1.} In \cite{AL}, an element $f\in \H$ was said to be of logarithmic
ergodic rate $\alpha \ge 0$, if $\sup_{n\ge 1}\frac{(\log n)^{1+\alpha}}{n}
\|U_n(f)\|<+\infty$. It was shown in \cite{AL} that if $f$ has logarithmic
ergodic rate $\alpha>0$ (and $T$ is unitary) then $(i)$ and $(ii)$ are valid.\\
{\bf 2.} Gaposhkin asked (\cite{Gap3}, P. 254) whether, in the case where
$T$ is induced by a measure
preserving transformation on $L^2(X,\mu)$, the convergence in norm in $(i)$
implies the almost sure convergence.
Point $(iii)$ may help to solve that question. In particular, one can deduce
from $(iii)$ that it sufficies to prove the a.s. convergence of $\{\sum_{k=1}
^{2^{2^n}}\frac{T^kf}{k}\}$. Indeed we have
$$
\max_{2^{2^n}+1\le k \le 2^{2^n+1}} |\sum_{l=2^{2^n}+1}^{k}
\frac{T^lf}{l}|\le \sum_{l=2^{2^n}+1}^{2^{2^{n+1}}} \frac{|U_lf|}{l^2}+
\max_{2^{2^n}+1\le k \le 2^{2^n+1}}\frac{|U_k|}{k},
$$
where the last term converges to 0 a.s. by the ergodic theorem
(or by \cite{Gap3}, p. 281). On the other hand, we have, by Cauchy-Schwarz
$$
Z_n^2:= (\sum_{l=2^{2^n}+1}^{2^{2^{n+1}}} \frac{|U_lf|}{l^2})^2
\le \big(\sum_{l=2^{2^n}+1}^{2^{2^{n+1}}}\frac{1}{k\log k} \big)\big(
\sum_{l=2^{2^n}+1}^{2^{2^{n+1}}} \frac{\log k |U_k|^2}{k^3}\big),
$$
which yields the result by $(iii)$ and Beppo-Levi's theorem.

\section{Behaviour of certain power series on the unit disk at the point $1$}

In this section we give the behaviour at $1$ of certain power series on the
unit disk when we know the behaviour of the coefficients. We did not find
these results in the litterature where people usually consider the case of the
circle or the interval $]-1,1[$. However we will closely follow a proof of
Zygmund in the case of the circle. We will state the results here and leave
the proofs to the appendix.\\

Following Zygmund, we say that $b~:~]u_0,+\infty[\rightarrow ]0,+\infty[$ is
slowly varying if for any $\delta >0$, $u^\delta b(u)$ and $u^{-\delta}b(u)$
are respectively non decreasing, non increasing, for $u$ large enough.
In particular, for every $k>0$, $b(ku)\underset{+\infty}\sim b(u)$
( even uniformly in $k\in [\eta, 1/\eta]$, $0<\eta <1$). Hence our
slowly varying functions are a particular case of the usual concept of slowly
varying functions.\\

Following Zhao-Woodroofe \cite{ZW}, define for any slowly varying function $b$
$$
\gamma_n:= \frac{c}{n}\sum_{k\ge n} \frac{1}{\sqrt{k^3b(k)}}, \qquad n\ge 1,
$$
where $c$ is chosen such that $\sum_{n\ge 1}\gamma_n=1$,
$$
B(z):= \sum_{n\ge 1}\gamma_n z^n \qquad \forall z\in \overline D.
$$
Then $B$ is well-defined and continuous, since the series is absolutely
converging. \\
Moreover $B$ is a convex combination of elements of $\overline D$, hence
$B(z)=1$ if and only if $z=1$, and $A:=\frac{1}{1-B}$ defines a continuous
function on $\overline D-\{1\}$ analytic on $D$. Hence there exists
$\{\alpha_n\}$ such that
\begin{equation}\label{A2}
A(z)=\frac{1}{1-B(z)}=\sum_{n\ge 0}\alpha_n z^n \qquad \forall z\in D.
\end{equation}

We will show in the appendix that
\begin{prop}\label{A1}
Let $b$ be any monotonic differentiable slowly varying function and let $A$
 be as above. We have
\begin{itemize}
\item [(i)]  $\displaystyle |A(z)|\underset{1}\sim
\frac{\sqrt{b(\frac{1}{|1-z|})}}{2c\sqrt \pi \sqrt{|1-z|}}$.
\item [(ii)] There exists $K>0$ such that $\sup_{n\ge 0}|\sum_{k=0}^n
\alpha_k z^k|\le K \frac{\sqrt{b(\frac{1}{|1-z|})}}{\sqrt{|1-z|}}$.
\item [(iii)] The series in \eqref{A2} converges on $\overline D-\{1\}$, and
the identity \eqref{A2} holds on $\overline D-\{1\}$.
\end{itemize}
\end{prop}

The proof uses the following proposition (of independent interest).

\begin{prop}\label{prop1}
Let $0<\beta< 1$ and $b$ be a slowly varying function. Then
$$
\sum_{n\ge 1} \frac{b(n)z^n}{n^\beta}= \Gamma (1-\beta)(1-z)^{\beta -1}
b(\frac{1}{|1-z|})+o(|1-z|^{\beta -1}b(\frac{1}{|1-z|})) \qquad
z\rightarrow 1,$$
where the power series is convergent on $\overline D-\{1\}$.
\end{prop}

The proof of  Proposition \ref{prop1} will be done in the appendix too.\\

\begin{theo}\label{norm}
Let $b$ be any  slowly varying function,  and let
 $\{\alpha_n\}$ be defined as above. For any normal contraction or
isometry $T$ on a Hilbert space $\H$ and any $f\in \H$ with spectral measure
$\mu_f$, the following are equivalent
\begin{itemize}
\item [(i)]   $\sum_{n\ge 0} \alpha_n T^n(f)$ converges in norm.
\item [(ii)]  $\displaystyle \int_{\overline D}\frac{b(\frac{1}{|1-z|})}
{|1-z|}\mu_f(dz)<+\infty$.
\item [(iii)]  $\sum_{n\ge 1} \frac{b(n)\|U_n(f)\|^2}{n^2} <+\infty$.
\end{itemize}
\end{theo}
\noindent {\bf Proof:}\\
Let $T$ be a normal contraction of a Hilbert space $\H$ and $f\in \H$.\\
By Proposition \ref{A1}, we can apply Proposition
\ref{CoL} with $\psi(x)=\frac{\sqrt{xb(x)}}{2c\sqrt \pi}$, to show that
$(i)$ and $(ii)$ are equivalent.\\
Hence it remains to show that $(ii)$ is equivalent to $(iii)$.
We will apply Lemma \ref{lem} with $\chi(x):= xb(x)$.\\
We just need to show that there exsits $\tau>1$, such that
$\chi(2^{n+1})\ge \tau \chi(2^n)$, for every $n\ge 1$. We have
$$
\chi(2^{n+1})=2\chi(2^n) \frac{b(2^n)}{b(2^{n+1})},
$$
and, since $b$ is slowly varying, $\frac{b(2^n)}{b(2^{n+1})}
\underset{n\rightarrow +\infty}\rightarrow 1$.\\
So the theorem is proved for normal $T$. In case $T$ is an isometry, the
result follows by the unitary dilation as in Lemma 2.3 of
\cite{CL}. \hfill $\square$

\noindent {\bf Remarks:} \\
{\bf 1.} It was proved in \cite{ZW} (see the proof of
Proposition 4), that a sufficient condition for $(i)$ is $\sum_{n\ge 1}
\frac{\sqrt{b(n)}\|S_n(f)\|}{n^{3/2}}<+\infty$.\\
{\bf 2.} Using Proposition 2.2, we can characterize the fact that
$f\in (I-T)^\alpha \H$, $0<\alpha<1$ (see \cite{DL} for the definition)
by $\sum_{n\ge 1}\frac{\|S_n(f)\|^2}{n^{3-2\alpha}}$.

\section{Ergodic theorems with rate}

We now give some ergodic theorems with rates as applications of the previous
section and of a result from \cite{ZW}, inspired by a method of \cite{DL}.
A different extension of \cite{DL} may be found in \cite{CL}.

Let us recall first the result of Zhao and Woodroofe \cite{ZW}.
For a contraction $T$ of a Hilbert space $\H$, define
$$
A(T)(f):=\sum_{n\ge 0}\alpha_nT^n(f),
$$
whenever the series converging, where $\{\alpha_n\}$ is defined by
\eqref{A2}.
\begin{theo}[Zhao-Woodroofe, \cite{ZW}]\label{ZW}
Let $b$ be any non decreasing slowly varying function. Let $\theta$ be a
measure preserving transformation of $(X,\Sigma,\nu)$ and $T$ be the
isometry induced by $\theta$ on $L^2(X,\nu)$. Let $f\in L^2(X,\nu)$ such
 that the series giving $A(T)(f)$ converges in  $L^2(X,\nu)$. Then
\begin{equation}\label{erg3}
\frac{1}{\sqrt {nb^*(n)}}\sum_{k=1}^nT^kf \underset
{n\rightarrow +\infty}\rightarrow 0 \qquad \mbox{$\nu$-a.s.},
\end{equation}
where $b^*(n):=\sum_{k=1}^n\frac{1}{kb(k)}$.
\end{theo}
\noindent {\bf Remarks:}\\
{\bf 1.} The theorem, as it is stated does not appear in \cite{ZW}, but it
is an easy
consequence of their Proposition 4 and Theorem 2, see also the proofs
therein.\\
{\bf 2.} It can be checked that the conclusion of Theorem \ref{ZW}
holds for slowly varying function (with our definition) without the non
decreasing assumption, under the extra assumption
$b(n)b^*(n)\underset{n\rightarrow +\infty}\rightarrow +\infty$, which is
 needed to ensure (18) in \cite{ZW}.
In particular, \eqref{erg3} remains valid for $b$ of the form $b(x)=
\frac{1}{(\log x)^\alpha (\log \log x)^\beta}$.\\
{\bf 3.} As in \cite{DL}, the proof of Theorem \ref{ZW} applies for
$T$ a Dunford-Schwarz operators (that is $T$ a contraction of every
$L^p(X,\Sigma,\nu)$, $(1\le p\le \infty)$.\\

Using Theorem \ref{norm} we deduce

\begin{theo}\label{ergo1}
Let $\theta$ be a measure preserving transformation of $(X,\Sigma,\nu)$  and
$T$ be the operator induced by $\theta$ on $L^2(X,\Sigma,\nu)$.
Let $b$ any non decreasing  slowly varying function. Let $f\in L^2(X)$,
such that
$$
\sum_{n\ge 1}\frac{b(n)\|\sum_{k=1}^nT^kf\|_{L^2(\nu)}^2}
{n^2}<+\infty.
$$
Then
$$
\frac{1}{\sqrt {nb^*(n)}}\sum_{k=1}^nT^kf \underset
{n\rightarrow +\infty}\rightarrow 0 \qquad \mbox{$\nu$-a.s.},
$$
where $b^*(n):=\sum_{k=1}^n\frac{1}{kb(k)}$.
\end{theo}
\noindent {\bf Proof:}\\
Apply Theorem \ref{norm} and Theorem \ref{ZW}.\\
{\bf Remarks:} \\
{\bf 1.} The theorem is also true with $b$ of
the form $b(x)=\frac{1}{(\log x)^\alpha (\log \log x)^\beta}$, see the
previous remark.\\
{\bf 2.} By the previous remark 3, the theorem is valid for $T$ a
Dunford-Schwarz operator, such that the restriction of $T$ to
$L^2(X,\Sigma, \nu)$ is either an isometry or a normal contraction (for
instance a normal Markov operator).

\begin{theo}\label{ergo2}
Let $\theta$ be a measure preserving transformation of $(X,\Sigma,\nu)$.
Let $f\in L^2(\nu)$ and  $\delta>1$, such that
\begin{equation}\label{erg1}
\sum_{n\ge 3}\frac{\log n (\log \log n)^\delta
\|\sum_{k=1}^nf\circ \theta^k\|_{L^2(\nu)}^2}{n^2}<+\infty.
\end{equation}
Then
$$
\frac{1}{\sqrt {n}}\sum_{k=1}^nf\circ \theta^k\underset
{n\rightarrow +\infty}\rightarrow 0 \qquad \mbox{$\nu$-a.s.}
$$
and the series
$$
\sum_{n\ge 1} \frac{f \circ \theta^n}{\sqrt n}
$$
converges $\nu$-a.s.
\end{theo}

\noindent {\bf Proof:}\\
Apply the previous theorem with $b=\log  (\log \log )^\delta$. Then
$\sum_{n\ge 1}\frac{1}{b(n)}<+\infty$ and the first assertion follows.
To prove the second one, notice that
\begin{equation*}\label{dec}
\sum_{k=1}^n \frac{f \circ \theta^k}{\sqrt k}=\sum_{k=1}^{n-1}U_k(\frac{1}
{\sqrt k}-\frac{1}{\sqrt{k+1}})+\frac{U_n}{\sqrt n}.
\end{equation*}
The second term is converging to $0$ by the previous result and the series on
the right hand side is $\nu$-a.s. absolutely converging since,
by Cauchy-Schwarz,
\begin{gather*}
\sum_{n\ge 1}|U_n||\frac{1}
{\sqrt n}-\frac{1}{\sqrt{n+1}}|\le C\sum_{n\ge 1}\frac{1}{\sqrt {n\log n
(\log \log n)^\delta}}\frac{|U_n|\sqrt {n\log n(\log \log n)^\delta}}
{n}\\
\le (\sum_{n\ge 1}\frac{1}{n\log n(\log\log n)^\delta})^{1/2}
(\sum_{n\ge 1}\frac{\log n (\log \log n)^\delta |U_n|^2}{n^2})^{1/2},
\end{gather*}
and the second series is converging by Beppo-Levi's theorem.
\hfill $\square$\\
{\bf Remarks.}\\
{\bf 1.} We could give a more general result by taking any
slowly varying $b$ such that $\sum_{n\ge 1}\frac{1}{nb(n)}<+\infty$.\\
{\bf 2.} It was proved in \cite{CL} (Theorem 3.3) that the conclusion
of Theorem \ref{ergo1} holds under the condition $\sup_{n\ge 3}
\frac{(\log n)^{3/2} (\log \log n)^\tau }{\sqrt n}\|\sum_{k=1}^nf\circ
\theta^k\|_{L^2(\nu)}<+\infty$, for some $\tau>1$.\\
{\bf 3.} In \cite{ZW}, the condition $\sum_{n\ge 3}
\frac{(\log n)^{1/2} (\log \log n)^\delta}{n^{3/2}}\|\sum_{k=1}^nf\circ
\theta^k\|_{L^2(\nu)}<+\infty$, for some $\delta>1/2$ was shown to be
sufficient.\\
{\bf 4.} Our condition is better than the ones from
\cite{CL} and \cite{ZW}, and those two seems to be uncomparable.
To see that our condition is better than the one in \cite{ZW}, using
the subadditivity of $\{\|U_n(f)\|_2\}$ as in \cite{PU}, one obtain
$\frac{n}{8}\frac{(\log n)^{1/2}\|U_nf\|}{n^{3/2}}\le
\sum_{n/4\le k\le n/2}\frac{(\log n)^{1/2}(\|U_kf\|+\|U_{n-k}(f)\|}{n^{3/2}}
\le \sum_{n/4\le k\le 3n/4}\frac{(\log k)^{1/2}\|U_kf\|}{k^{3/2}}
\underset{n\rightarrow +\infty }\rightarrow 0$,
which, combined again with the condition in {\bf 3.} gives \eqref{erg1}.\\
{\bf 5.} A good condition was obtained by Wu in  \cite{Wu}, based on the
dyadic chaining. His proof leads to a weak type maximal inequality
and works in $L^p$ spaces. It was shown there that for all $\eta>0$,
$$
\displaystyle \sum_{n\ge 1}\nu\{ \max_{k=1}^{2^n}|\sum_{l=1}^kf\circ
\theta^k|\ge
\eta \sqrt {2^n}\}\le \frac{2}{\eta^2}\bigg(\sum_{j\ge 0}\bigg(
\frac{\|\sum_{k=1}^{2^j}f\circ \theta^k\|_{L^2(\nu)}^2}{2^j}
\bigg)^{1/3}\bigg)^3.
$$
We were not able to compare our condition to the one of \cite{Wu}
but in the applications we will see that our condition yields better results.

\begin{theo}\label{ergo3}
Let $\theta$ be a measure preserving transformation of $(X,\Sigma,\nu)$.
Let $f\in L^2(\nu)$ such that
$$
\sum_{n\ge 3}\frac{\log n
\|\sum_{k=1}^nf\circ \theta^k\|_{L^2(\nu)}^2}{n^2}<+\infty.
$$
Then
$$
\frac{1}{\sqrt {n \log \log n}}\sum_{k=1}^nf\circ \theta^k\underset
{n\rightarrow +\infty}\rightarrow 0 \qquad \mbox{$\nu$-a.s.}
$$
\end{theo}
\noindent {\bf Proof:}\\
Apply Theorem \ref{ergo1} with $b=\log  $. Then
$\sum_{k= 2}^n \frac{1}{n\log n} \underset{n\rightarrow +\infty}\sim
\log \log n$. \hfill $\square$\\

Finally, we obtain

\begin{theo}\label{ergo4}
Let $\theta$ be a measure preserving transformation of $(X,\Sigma,\nu)$.
Let $f\in L^2(\nu)$ such that there exists $\beta>1$ such that
$$
\sum_{n\ge 3}\frac{
\|\sum_{k=1}^nf\circ \theta^k\|_{L^2(\nu)}^2}{n^2\log n(\log \log n)^\beta}
<+\infty.
$$
Then
\begin{equation}\label{erg}
\frac{1}{\sqrt {n} \log n(\log \log n)^{\beta/2}}\sum_{k=1}^nf\circ \theta^k
\underset {n\rightarrow +\infty}\rightarrow 0 \qquad \mbox{$\nu$-a.s.}
\end{equation}
\end{theo}
\noindent {\bf Proof:}\\
Apply Theorem \ref{ergo1}, see the remark after the theorem, with $b(x)=
\frac{1}{\log x (\log \log x)^\beta}$. Then $b^*(n)=\sum_{k=1}^n
\frac{\log n (\log \log n)^\beta}{n}\sim (\log n)^2 (\log \log n)^\beta$.\\
{\bf Remarks:}\\
{\bf 1.} It follows from Theorem \ref{ergo4}, that \eqref{erg} holds for
every $\beta>1$ for every $f$ satisfying
\begin{equation}\label{sup}
\sup_{n\ge 1}\frac{\|\sum_{k=1}^nf\circ \theta^k\|_{L^2(\nu)}}{\sqrt n}<
+\infty.
\end{equation}
{\bf 2.} Under \eqref{sup}, the conclusion \eqref{erg} was shown to hold
with a power $\alpha >3/2$ on the logarithm in \cite{CoL1} and \cite{Weber}
for Dunford-Schwarz operators and in \cite{Gap2} for unitary operators.

\section{Applications to the Quenched CLT and the LIL}

We now give some applications of the previous section to the study of
ergodic stationary processes. The main purpose is to obtain conditions
under
which the process may be approximated by a martingale such that
the remainder will satisfy the assumptions of Theorems
\ref{ergo2} or \ref{ergo3}. Then we can obtain limit theorems such
as the quenched CLT or Law of the Iterated Logarithm (LIL).
\medskip

Let $\{X_n\}$ be a stationary ergodic process in $L^2(\Omega,\F,\P)$, with natural
 filtration $\{\F_n\}$  and $\theta$ the shift associated. Write
$S_n=X_1+\ldots +X_n$.

\smallskip

\noindent {\bf Definition:} \emph{We say that $\{X_n\}$ admits a martingale
approximation
if there exists $M\in L^2$ such that  $\{M_n:=\sum_{k=1}^nM\circ
\theta^k\}$ is a martingale with respect to $\{\F_n\}$ and, if $R_n:=S_n-M_n$,
satisfies to $\|R_n\|=o(\sqrt n)$.}

\medskip

Clearly, if $\{X_n\}$ admits a martingale approximation, then $\{M_n\}$ is
unique and
\begin{equation}\label{sigma}
\sigma^2:=\E[M^2]=\lim_{n\rightarrow +\infty} \frac{\E[S_n^2]}{n}.
\end{equation}

In particular, the coming construction yields to the same martingale
approximation, but the way to obtain it gives different estimates.

We will be also concerned with the quenched central limit theorem for
Markov chains. \\

Let $\{W_n\}_{n\ge 0}$ be a stationary ergodic Markov chain
with state space $({\mathbb S},{\mathcal S})$, transition probability $P$,
invariant initial distribution $m$, and corresponding Markov operator $P$ on
$L^2({\mathcal S},m)$. For $x\in {\mathbb S}$, denote by $\P_x$ the probability of
the chain starting from $x$, defined on the product $\sigma$-algebra of
$\Omega:={\mathbb S}^{\bf N}$. We will write $\P_m$ the probability of the
chain
starting according to $m$.

Thoughout the paper the Markov chain will be as above and we will use
these notations.

Recall that, since $m$ is invariant, the chain $\{W_n\}_{n\ge 0}$ may be
extended to a chain indexed by ${\bf Z}$.

\noindent {\bf Definition:}  \emph{Let $\{W_n\}$ be a Markov chain as above,
and $f\in L^2( {\mathcal S},m)$. We say that $\{f(W_n)\}$
 satisfies the} quenched CLT and invariance principle, \emph{if, for
$m$-almost every $x\in {\mathcal S}$, the sequence $\{
\frac{1}{\sqrt n}\sum_{k=1}^nf(W_k)\}$ converges in distribution, in the space
$(\Omega,\P_x)$ to a (possibly degenerate) Gaussian distribution
$\N(0,\sigma(f)^2)$ (with variance $\sigma(f)^2$ independent of $x$), and if
also  the invariance principle holds.}
\\

\subsection{Using Wu's estimates for the martingale approximation}

Let $\{X_n\}_{n\in {\bf Z}}$ be an ergodic  stationnary process on a probability space
$(\Omega, \F,\P)$ with $\|X_0\|_2<+\infty$. Let  $\{\F_n\}_{n\in {\bf Z}}$ be
the natural filtration of $\{X_n\}$ and $\theta$ the shift associated with
$\{X_n\}$.\\

Following Wu \cite{Wu}, we define
\begin{itemize}
\item [$\bullet$] $S_n:=\sum_{k=1}^n X_k \qquad \forall n\ge 1$.
\item [$\bullet$] $Z_{i,k}:= \E[X_i|\F_k]-\E[X_i|\F_{k-1}] \qquad \forall
i\ge k $.
\item [$\bullet$] $\Theta_m:=\sum_{i\ge m} \|Z_{i,0}\|_2 \qquad
\forall m\ge 0$.
\item [$\bullet$] $D_k:= \sum_{i\ge k}Z_{i,k} \qquad \forall k\ge 1$
\quad (whenever the series converges $\P$-a.s.)
\item [$\bullet$] $M_n:=\sum_{k=1}^nD_k \qquad \forall n\ge 1$.
\item [$\bullet$] $R_n:=S_n-M_n \qquad \forall n\ge 1$.
\end{itemize}
Clearly, if $\Theta_0<+\infty$, then $\{D_k\}$ is well defined and, for every
$k \ge 0$, $D_k=D_0\circ \theta ^k$.
Moreover $\{M_n\}$ is a stationary martingale (with ergodic increments).
Wu proved the following
\begin{theo}[Wu, \cite{Wu}]\label{Wu}
Let $\{X_n\}$ be an ergodic  stationary process with $\E [X_0^2]<+\infty$
 and
$\E[X_0]=0$. If $\Theta_0<+\infty$, then there exists $K>0$ such that
$$
\E[R_n^2]\le K \sum_{j=1}^n \Theta_j^2.
$$
\end{theo}
\noindent {\bf Remark.} The theorem in \cite{Wu} is stated when $\{X_n\}$ is
the functional of a Markov chain but the proof uses only the stationarity
of the process.

We deduce the following
\begin{theo}\label{rate}
Let $b$ be any  slowly varying function. Let $\{X_n\}$ be an ergodic
 centered
stationary process in $L^2(\Omega,\F,\P)$ such that
$$
\sum_{n\ge 1}\frac{b(n) \Theta_n^2}{n}<+\infty.
$$
Then
$$\displaystyle \frac{R_n}{\sqrt {nb^*(n)}} \underset{n\rightarrow +\infty}
\rightarrow 0\qquad \mbox{$\P$-a.s.},
$$
where $b^*(n):=\sum_{k=1}^n \frac{1}{kb(k)}$.
\end{theo}
\noindent {\bf Proof:}\\
Define $V_1=R_1$ and for every $n\ge 2$, $V_n:=R_n-R_{n-1}$.
Then $V_{n+1}=V_1\circ \theta^n$ and $R_n=\sum_{k=0}^{n-1} V_1\circ \theta^k$.
Hence, by Theorem \ref{ergo1} and Theorem \ref{Wu},
the conclusion of theorem \ref{rate} holds as soon as
$$
\sum_{n\ge 1} \frac{b(n)}{n^2}\sum_{k=1}^n\Theta_k^2 <+\infty.
$$
But, using that $b(n)/\sqrt n$ is non increasing at infinity, we have
\begin{gather*}
\sum_{n\ge 1} \frac{b(n)}{n^2}\sum_{k=1}^n\Theta_k^2 =
\sum_{k\ge 1} \Theta_k^2 \sum_{n\ge k}\frac{b(n)}{\sqrt n}\frac{1}{ n^{3/2}}
\le C \sum_{k\ge 1} \Theta_k^2 \frac{b(k)}{ k},
\end{gather*}
which proves the theorem. \hfill $\square$

\begin{theo}\label{LIL}
Let $\{X_n\}$ be an ergodic centered
stationary process in $L^2(\Omega,\F,\P)$ such that
\begin{equation}\label{WC}
\sum_{n\ge 1}\frac{\log n \Theta_n^2}{n}<+\infty.
\end{equation}
Then
$$
\limsup_{n\rightarrow +\infty}\frac{S_n}{\sqrt{2n\log \log n}}=\sigma \qquad
\mbox{$\P$-a.s.},
$$
where $\sigma$ is defined in \eqref{sigma}
\end{theo}
\noindent {\bf Proof:} \\
By the previous theorem with $b:=\log$,
$$
\limsup_{n\rightarrow +\infty}\frac{R_n}{\sqrt{n\log \log n}}=0 \qquad
\mbox{$\P$-a.s.}
$$
The result follows then from the identity $S_n=M_n+R_n$ and Stout's law of the
iterated logarithm (see \cite{Stout}) applied to the martingale $\{M_n\}$.
 \hfill $\square$
\begin{prop}\label{propC}
Let $\{X_n\}$ be an ergodic centered stationary process in
$L^2(\Omega,\F,\P)$, such
that $\sum_{n\ge 1}(\log n)^3 \|\E[X_n|\F_0]\|^2<+\infty$, then \eqref{WC}
holds.
\end{prop}
\noindent {\bf Proof:}\\
Using Cauchy-Schwarz, we have
\begin{gather}
\nonumber\Theta_n^2=\big(\sum_{k\ge n}(\|\E[X_k|\F_0]\|^2-\|\E[X_{k+1}|
\F_0]\|^2)^{1/2}\big)^2\\
\nonumber \le (\sum_{k\ge n}\frac{1}{k(\log k)^2})(\sum_{k\ge n} k(\log k)^2
(\|\E[X_k|\F_0]\|^2-\|\E[X_{k+1}|\F_0]\|^2))\\
\label{test1}\le \frac{C_1}{\log n} \sum_{k\ge n} k(\log k)^2
(\|\E[X_k|\F_0]\|^2-\|\E[X_{k+1}|\F_0]\|^2).
\end{gather}
On the other hand, using Abel summation
\begin{gather}
\nonumber \sum_{k\ge n} k(\log k)^2 (\|\E[X_k|\F_0]\|^2-
\|\E[X_{k+1}|\F_0]\|^2)\\
\nonumber = n(\log n)^2\|\E[X_n|\F_0]\|^2+ \sum_{k\ge n+1}\|\E[X_k|\F_0]\|^2
(k(\log k)^2-(k-1)(\log (k-1))^2)\\
\label{test2}\le n(\log n)^2\|\E[X_n|\F_0]\|^2+ C\sum_{k\ge n+1}
\|\E[X_k|\F_0]\|^2
(\log k)^2.
\end{gather}
Moreover,
$$
\sum_{n\ge 2}\frac{1}{n} \sum_{k\ge n+1}\|\E[X_k|\F_0]\|^2
(\log k)^2= \sum_{k\ge 3} \|\E[X_k|\F_0]\|^2(\log k)^2 \sum_{n=2}^{k-1}
\frac{1}{k} \le K \sum_{k\ge 3} \|\E[X_k|\F_0]\|^2(\log k)^3,
$$
which, combined with \eqref{test1} and \eqref{test2}, yields the desired
result. \hfill $\square$
\begin{theo}\label{quenched}
Let $\{W_n\}$ be a stationary ergodic Markov chain. Let $f\in
L^2({\mathcal S},m)$, such that  there exists $\delta>1$ such that
$$
\sum_{n\ge 1} \frac{\log n (\log \log n)^\delta}{n}
\big(\sum_{k\ge n}
(\|P^{k-1}f\|_{L^2(m)}^2-\|P^{k}f\|_{L^2(m)}^2)^{1/2}\big)^2<+\infty,
$$
then $\{f(W_n)\}$ statisfies the quenched CLT and invariance principle.
\end{theo}
\noindent {\bf Proof:}\\
By Theorem \ref{rate} with $b:=\log (\log \log)^\delta$, we have
$$
\limsup_{n\rightarrow +\infty}\frac{R_n}{\sqrt{n}}=0 \qquad
\mbox{$\P$-a.s.},
$$
where $S_n=M_n+R_n$ and $\{M_n\}$ is a stationary martingale.
The end of the proof is now similar to  \cite[p. 75]{DL2} . \hfill $\square$\\
{\bf Remark.}\\
The conclusion of the theorem may be deduced from Wu \cite{Wu} under the
condition $\sum_{n\ge 1}\frac{\Theta_n^{2/3}}{n}$. Since $\{\Theta_n\}$ is non
increasing, that condition implies that $\Theta_n=O((\log n)^{-3/2})$, which
implies that our condition is satisfied.\\

\subsection{Using Kipnis-Varadhan's method for martingale approximation}

Kipnis and Varadhan discovered a useful way to obtain martingale
approximation in the case of a Markov chain. They worked with Markov chains
with symmetric Markov operator, but the method applies in general,
as it was done by Maxwell-Woodroofe \cite{MW}. It is mentionned in \cite{MW}
 that such an approximation may be obtained (from the Markov chain case)
 for a general stationary process
$\{X_n\}$ by considering the process $\{W_n\}$ defined as
$W_n:= (...,X_{n-1},X_n)$. We will show that actually  the method applies
directly to general stationary processes. See also (\cite{MPU}, p. 9) for a
general (related) scheme of approximation by martingale. \\

Let $\{X_n\}_{n\in {\bf Z}}$ be an ergodic  stationary process in $L^2(\Omega,\F,\P)$.
let $S_n:=X_1+\ldots X_n$. Define
$$
\Gamma_t:= \sum_{n\ge 0}t^nX_n \qquad \forall 0\le t<1,
$$
where the series converges  in $L^2$, for every $t\in [0,1[$. Moreover, taking
$\{t_k\}\subset [0,1)$ converging to $1$ and using Beppo-Levi's theorem, there
 exists $\Omega_0\in \F$, with $\P(\Omega_0)=1$, such that for every
$\omega \in \Omega_0$, the series giving $\Gamma_t(\omega)$ is absolutely
convergent for every $t\in [0,1)$.\\
Then
\begin{equation}\label{X_0}
\E[\Gamma_t|\F_0]=\sum_{n\ge 0}t^n\E[X_n|\F_0] \mbox{ and }
X_0=\E[\Gamma_t|\F_0]-t\E[\Gamma_t\circ \theta |\F_0].
\end{equation}
Putting $t^n=(1-t)\sum_{k\ge n}t^k$ in $\Gamma_t$ and using Fubini, we
obtain $\Gamma_t=(1-t)\sum_{n\ge 0}t^nS_n$. Hence
\begin{equation}\label{Abel}
\|\E[\Gamma_t|\F_0]\|_2\le (1-t)\sum_{n\ge 0}t^n \|\E[S_n|\F_0]\|_2.
\end{equation}
Let $\theta$ be the shift associated with $\{X_n\}$ and
$\{\F_n\}$ be its natural filtration. Define
$$
\varphi_t:=(\E[\Gamma_t|\F_0])\circ \theta-\E[\Gamma_t\circ \theta|\F_0]
=(\E[\Gamma_t\circ \theta |\F_1])-\E[\Gamma_t\circ \theta|\F_0]
$$
Then $M_n(t):=\sum_{k=1}^n \varphi_t\circ \theta^k$ is a martingale with
(ergodic) stationary increments. \\

We would like to show that under some estimates on
$\{\|\E[S_n|\F_0]\|_2\}$, $\{M_n(t)\}$ converges in $L^2$, when
$t$ goes to $1$, to a martingale $\{M_n\}$ with stationary increments such
that one can estimate $\{\|S_n-M_n\|_2\}$.

\begin{prop}\label{prop}
Let $\{X_n\}_{n\in {\bf Z}}$ be an ergodic  stationary process such that
$$
\sup_{n\ge 1} \frac{(\log n)^2 (\log \log n)^\tau}{\sqrt n} \|\E[S_n|\F_0]\|_2
<+\infty,
$$
for some $\tau \in {\bf R}$. \\
Then $\{X_n\}$ admits a martingale approximation $\{M_n\}$ such that
$$
\sup_{n\ge 1} \frac{\log n (\log \log n)^\tau}{\sqrt n} \|S_n-M_n\|_2
<+\infty.
$$
\end{prop}
\noindent {\bf Proof:} see the appendix.

\begin{theo}\label{lil}
Let $\{X_n\}$ be an ergodic  centered
stationary process in $L^2(\Omega,\F,\P)$ such that there exists $\tau>1/2$
\begin{equation}\label{ZWC}
\sup_{n\ge 1}\frac{(\log n)^2 (\log \log n)^\tau}{\sqrt n}
\|\E[S_n|\F_0]\|_2<+\infty.
\end{equation}
Then
$$
\limsup_{n\rightarrow +\infty}\frac{S_n}{\sqrt{2n\log \log n}}=\sigma \qquad
\mbox{$\P$-a.s.},
$$
where $\sigma$ is defined by \eqref{sigma}.
\end{theo}
\noindent {\bf Proof:}\\
By Proposition \ref{prop}, there exists a martingale $\{M_n\}$ and $R\in L^2
(\Omega,\P)$, such that, $R_n:=S_n-M_n=\sum_{k=1}^n R\circ \theta^k$ and
$\sup_{n\ge 1} \frac{\log n (\log \log n)^\tau}{\sqrt n} \|R_n\|_2
<+\infty$. Then apply, Theorem \ref{ergo3} with $f:=R$ to obtain
$\frac{R_n}{\sqrt{n\log\log n}}\underset{n\rightarrow}\rightarrow 0$.
Thus the result follows from Stout's LIL for martingale, see \cite{Stout}.
\hfill $\square$\\
{\bf Remark:}\\
Zhao and Woodroofe \cite{ZW} obtained the conclusion of the theorem under the
condition $\sum_{n\ge 1} \frac{(\log n)^{3/2}}{n\sqrt n} \|\E[S_n|\F_0]
\|_2<+\infty$. Using the subadditivity $\{\|\E[S_n|\F_0]\|_2\}$, it can be
shown that the previous condition implies $\frac{(\log n)^{3/2}}
{\sqrt n} \|\E[S_n|\F_0]\|_2\underset{n\rightarrow +\infty}\rightarrow 0$.
If one knew that $\{\frac{\|\E[S_n|\F_0]\|_2}{\sqrt n}\}$ is non increasing,
one could deduce from the Zhao-Woodroofe condition that $\frac{(\log n)^{5/2}}
{\sqrt n} \|\E[S_n|\F_0]\|_2\underset{n\rightarrow +\infty}\rightarrow 0$,
which implies our condition.
Of course there is no reason why this monotony assumption should be
satisfied but from a practical point of view one has to estimate
$\|\E[S_n|\F_0]\|_2$ and such an estimation often leeds to monotony.


\begin{theo}\label{quench}
Let $\{W_n\}$ be a stationary ergodic Markov chain. Let $f\in
L^2({\mathcal S}),m)$ such that there exists $\delta>1$ such that
$$
\sup_{n\ge 1} \frac{(\log n)^2 (\log \log n)^\delta}{\sqrt n}
\|\sum_{k=1}^nP^kf\|_2 <+\infty
$$
then $\{f(W_n)\}$ satifies the quenched CLT and invariance principle.
\end{theo}
\noindent {\bf Proof:}\\
Apply Proposition \ref{prop} to the process $\{f(W_n)\}$, noticing that
$\E[f(W_n)|\F_0]=P^nf(W_0)$, to obtain a martingale $\{M_n\}$ such
that $\sup_{n\ge 1} \frac{\log n (\log \log n)^\delta}{\sqrt n}
\|S_n-M_n\|_2 <+\infty$. Then apply Theorem \ref{ergo2} and the proof may be
finished as in \cite[p. 76]{DL2}. \hfill $\square$\\
{\bf Remarks:}\\
This theorem improves a theorem of \cite{CL} where we obtained a  condition
with a power $5/2$ on the logarithm. The condition
$\sum_{n\ge 1} \frac{(\log n)^{3/2}(\log \log n)^\tau}{\sqrt n}
\|\sum_{k=1}^nP^kf\|_2<+\infty$ for some $\tau>1$ was shown to be
sufficient in \cite{ZW}.

\section{The case of Markov chains with normal transition operator}

Let $\{W_n\}_{n\ge 0}$ be a stationary ergodic Markov chain with state space
$({\mathbb S},{\mathcal S})$, transition probability $P$, invariant initial
distribution $m$, and corresponding Markov operator $P$ on
$L^2({\mathcal S},m)$. For $f\in L^2({\mathcal S},m)$, denote $S_n(f):=
\sum_{k=1}^nf(W_k)$.\\

Throughout this section we assume that $P$ is a normal operator
on $L^2({\mathcal S},m)$, that is $PP^*=P^*P$. \\
Define  $U_n(f):=\sum_{k=1}^n P^k(f)$.
\begin{prop} \label{normal}
Let $\{W_n\}$ be a Markov chain with normal transitin operator as above
and $f\in L^2({\mathcal S},m)$, such that
$\sum_{n\ge 1}\frac{\|U_n(f)\|_2^2}{n^2}<+\infty$. Then $\{f(W_n)\}$ admits a
martingale approximation $\{M_n\}$ such that
\begin{equation}\label{estimate}
\frac{1}{n}\|S_n(f)-M_n\|_2^2\le C
\big(\frac{1}{n}\sum_{k=1}^n \frac{\|U_k(f)\|_2^2}{k}+\sum_{k\ge n+1}
\frac{\|U_k(f)\|_2^2}{k^2}\big).
\end{equation}
\end{prop}
\noindent {\bf Remark:} The fact that  $\|S_n(f)-M_n\|_2
=o(\sqrt n)$ follows from our assumption and Kronecker's lemma.
In view of Proposition \ref{sqrt}, we obtain that
whenever $f\in \sqrt{I-P}L^2({\mathcal S},m)$, $\{f(W_n)\}$ admits a
martingale
approximation. In particular, we recover that for
a normal $P$, the condition $f\in \sqrt{I-P}L^2({\mathcal S},m)$ implies the
central limit theorem (see \cite{KV} for $P$ symmetric and  \cite{DL0} for
$P$ normal).\\
{\bf Proof:} see the appendix.

\begin{theo}\label{theonorm}
Let $\{W_n\}$ be a stationary ergodic Markov chain (as above) whose
 Markov operator $P$ is normal on $L^2({\mathcal S},m)$. Let $f\in
L^2({\mathcal S},m)$.
\begin{itemize}
\item [(i)]  If there exists $\delta>1$ such that
\begin{equation*}\label{norm1}
\sum_{n\ge 1} \frac{(\log n)^2 (\log \log n)^\delta\|\sum_{k=1}^nP^kf\|_2^2 }
{n^2}<+\infty
\end{equation*}
then $\{f(W_n)\}$ satisfies the quenched CLT and invariance principle.
\item [(ii)] If
$$
\sum_{n\ge 1} \frac{(\log n)^2 \|\sum_{k=1}^nP^kf\|_2^2 }{n^2}<+\infty,
$$
then
$$
\limsup_{n\rightarrow +\infty} \frac{\sum_{k=1}^n f(X_k)}
{\sqrt {2n\log \log n}}=\sigma \qquad \mbox{$\P_m$-a.s.}
$$
\end{itemize}
\end{theo}
\noindent {\bf Proof:}\\
By proposition \ref{normal}, there exists  a martingale $\{M_n\}$ such that
$R_n:=S_n(f)-M_n$ satisfies \eqref{estimate}, and $R_n=\sum_{k=0}^{n-1}
R\circ \theta^k$.\\
Let $\delta \ge 0$. we have
\begin{gather*}
\sum_{n\ge 1} \frac{\log n (\log\log n)^\delta}{n^2}\|R_n\|^2\le
\sum_{n\ge 1}\frac{\log n (\log \log n)^{\delta}}{n^2}\big(
\sum_{k=1}^n\frac{\|U_k(f)\|_2^2}{k}+n\sum_{k\ge n}\frac{\|U_k(f)\|_2^2}{k^2})
\\
=\sum_{k\ge 1}\frac{\|S_k(f)\|_2^2}{k}\sum_{n\ge k}
\frac{\log n (\log \log n)^{\delta}}{n^2}+\sum_{k\ge 1}\frac{\|U_k(f)\|_2^2}{k}
\sum_{n=1}^k \frac{\log n (\log \log n)^{\delta}}{n}\\
\le C \sum_{k\ge 1} \frac{(\log k)^2
(\log \log k)^{\delta}\|U_k(f)\|_2^2}{k^2}.
\end{gather*}
Taking $\delta=0$ (respectively $\delta >1$), one may apply  Theorem
\ref{ergo3} (resp. Theorem \ref{ergo2}) to obtain
$$
\frac{1}{\sqrt {n\log \log n}}W_n\underset
{n\rightarrow +\infty}\rightarrow 0 \qquad \mbox{$\P_m$-a.s.}
$$
respectively
$$
\frac{1}{\sqrt {n}}W_n\underset
{n\rightarrow +\infty}\rightarrow 0 \qquad \mbox{$\P_m$-a.s.}
$$
Then the proof may be finished as in Theorem \ref{LIL} or Theorem
\ref{quenched}. \hfill $\square$

\begin{cor}
Assume that $\sum_{n\ge 1} (\log n)^2 (\log \log n)^\delta \|P^nf\|^2
<+\infty$. Then, if $\delta >1$, $(i)$ of Theorem \ref{theonorm}
is valid, and if $\delta=0$, $(ii)$ is valid.
\end{cor}
\noindent {\bf Proof:}\\
Denote $\mu_f$ the spectral measure of $f$ associated with $P$. We have
\begin{gather*}
\|\sum_{k=1}^nP^kf\|^2=\int_{D_1}|\sum_{k=1}^nz^k|^2\mu_f(dz)
\le 2 \sum_{k=1}^n \sum_{l=1}^k\int_{D_1}|z|^{k+l}\mu_f(dz)\\
\le 2 \sum_{k=1}^n \sum_{l=1}^k \|P^{[(l+k)/2]}f\|^2
\le 4 \sum_{k=1}^n \sum_{l=[(k+1)/2]}^k \|P^lf\|^2.
\end{gather*}
Hence
\begin{gather*}
\sum_{n\ge 1} \frac{(\log n)^2 (\log \log n)^\delta\|\sum_{k=1}^nP^kf\|_2^2 }
{n^2}\\
\le  4\sum_{k\ge 1}\sum_{l=[(k+1)/2]}^k \|P^lf\|^2 \sum_{n\ge k}
\frac{(\log n)^2 (\log \log n)^\delta}{n^2} \\
\le C \sum_{k\ge 1} \frac{(\log k)^2 (\log \log k)^\delta}{k}
\sum_{l=[(k+1)/2]}^k \|P^lf\|^2\\
\le C'\sum_{l\ge 1}  \|P^lf\|^2 \sum_{k=l}^{2l+1}\frac{(\log k)^2
(\log \log k)^\delta}{k}
\le C''\sum_{l\ge 1}  \|P^lf\|^2 (\log l)^2 (\log \log l)^\delta
\end{gather*}

\section{Comparison of the conditions of the theorems and examples}
We applied the results of section 2 to two different methods (Wu \cite{Wu}
and Zhao-Woodroofe \cite{ZW}) to obtain
martingale approximations. As we already mentionned the results from
section 2 essentially yields to weaker conditions than the ones obtained
in \cite{Wu} and \cite{ZW}. Now, we would like to compare our own
conditions.
\\

Let $\{\varepsilon_n\}_{n\in {\bf Z}}$ be iid centered random variables with
$\E [|\varepsilon_1|^2]<\infty$. Let $\{a_n\}_{n\ge 0}$ such that
$\sum_{n\ge 0}a_n^2<+\infty$. Define the linear process
$X_n:=\sum_{k\ge 0}a_k\varepsilon_{n-k}$.

\begin{prop}
There exists a linear process $\{X_n\}$ satisfying \eqref{WC} but which does
not satisfy \eqref{ZWC}.
\end{prop}
\noindent {\bf Proof:}\\
Define $a_i:= \frac{1}{k^{9/4}}$ if there exists $k\ge 1$ such that
$i=2^k$ and $a_i=0$ otherwise. Let $\{X_n\}$ be the associated linear process.
\\
Let show first that there exists $C>0$ such that
$$
\sup_{n\ge 3}\frac{(\log n)^{7/4} }{\sqrt n}
 \|\E[S_n|\F_0]\|_2 \ge C.
$$
Let $k\ge 1$, we have
\begin{gather*}
\|\E[S_{2^k}^2]\|^2     = \sum_{j=0}^\infty (a_{j+1}+\ldots a_{j+2^k})^2
    \ge  \sum_{l\ge k} \sum_{j=2^l-2^k}^{2^l-1} a_{2^l}^2\\
    \ge  2^k \sum_{l\ge k}a_{2^l}^2
   \ge \frac{2^k}{k^{7/2} },
\end{gather*}
which proves the above claim result.\\
On the other hand, for every $k\ge 2$, we have $\Theta_{2^k}=\sum_{n\ge 2^k}
a_n=\sum_{l\ge k}\frac{1}{l^{9/4}}=O(\frac{1}{k^{5/4}})$. Hence
$$
\sum_{k\ge 2}k\Theta_{2^k}^2 \le C \sum_{k\ge 2}\frac{1}{k^{3/2}}<+\infty,
$$
which implies the convergence of the series in \eqref{WC} since
$\{\Theta_n\}$ is non increasing.
\\

We now recall an example that we considered in \cite{CL}, in order to
compare theorems 3.3 and 5.2 in case of a Markov chain with normal transition
operator.  \\

Take $\alpha:=2{\rm e}$. Let $R_\alpha$ be the rotation of the unit circle of
angle $\alpha$ and define $P=P_\alpha:=\frac{1}{4}(2I+R_\alpha+R_{-\alpha})$.
Then $P$ is a symmetric operator (i.e. $P=P^*$).

\begin{prop}
There exists $f\in L^2[0,1]$ such that
\begin{equation*}
\sum_{n\ge 1}\frac{(\log n)^2}{ n^2}\|\sum_{k=1}^n
P_\alpha^kf\|^2<+\infty ,
\end{equation*}
and, for every $\delta>1$,
$$
\sum_{n\ge 1}\frac{\log n(\log \log n)^\delta}{n}(\sum_{k\ge n}
(\|P_\alpha^kf\|_2^2
-\|P_\alpha ^{k+1}\|_2^2)^{1/2})^2=+\infty.
$$
In particular, Theorem \ref{theonorm} $(ii)$ applies while Theorem
\ref{quenched} does not
\end{prop}
\noindent {\bf Proof:}\\
Let $f\in L^2[0,1]$, with Fourier expansion $f(x)=
\sum_{n\in {\bf Z}} c_n {\rm e}^{2i\pi n x}$. Then, for every $0\le x\le 1$,
$$
P^kf(x)=\sum_{n\in {\bf Z}}c_n \cos^{2k}(\pi n\alpha){\rm e}^{2i\pi nx}.
$$
If $n=l!$, for some $l\ge 3$, define $c_{-n}=c_n=\frac{1}{n^{3/2}(\log n)^2
}$, and define $c_n=0$ otherwise.\\
We have
\begin{gather*}
\|P^{k+1}f\|_2^2-\|P^kf\|_2^2= 2 \sum_{n\ge 3}c_{n!}^2 \cos^{4k}(\pi n! \alpha)
(1-\cos^4(\pi n!\alpha)).
\end{gather*}
It follows from Lemma 5.4 of \cite{CL} and from the proof of Lemma 5.5,
that there exists $C_1>0$ such that, for every $n\ge 2\pi$ and every
$k\le n^2$,
$1-\cos (\pi n!\alpha)\ge \frac{C_1}{n^2}$ and $\cos^{4k}(\pi n! \alpha)\ge
C_1$.
Hence we have, for every $k\ge 4\pi^2$
\begin{gather*}
\|P^{k+1}f\|_2^2-\|P^kf\|_2^2\ge \sum_{n\ge \sqrt k} \frac{C_1^2}{n^5(\log n)^4
}\\
\ge  \frac{C_2}{k^2 (\log k)^4 } .
\end{gather*}
Hence
$$
\sum_{k\ge n}(\|P_\alpha^kf\|_2^2-\|P_\alpha ^{k+1}\|_2^2)^{1/2}
\ge \frac{C_3}{\log k },
$$
which yields to the desired conclusion.\\
It remains to prove  \eqref{eq}. \\
By Lemma 5.5 of \cite{CL}, there exists $K>0$ such that for every $m\ge 1$
\begin{gather*}
\|\sum_{k=1}^m P_\alpha^k\|_2^2 \le K + K\sum_{7\le n\le \sqrt m}n^4
c_{n!}^2+Km^2\sum_{n>\sqrt m}c_{n!}^2.
\end{gather*}
Hence
\begin{gather*}
\|\sum_{k=1}^n P_\alpha^k\|_2^2 \le C+\frac{Cn}{(\log n)^4 },
\end{gather*}
which proves that \eqref{eq} is satisfied.\\

We now look at the case of $\rho$-mixing processes.
Let $\{X_n\}_{n\in {\bf Z}}$ be  a stationary process,
define $\rho(n):= \sup\{\|E[Y|\F_0]\|_2/\|Y\|_2~:~Y\in L^2(\F^n),~\E[Y]=0\}$,
where $\F_n=\sigma \{X_k,~k\ge n\}$. Then, by \cite[p. 15]{MPU}, we have
$$
\|\E[S_{2^{r+1}}|\F_0]\|_2\le C\sum_{j=0}^r 2^{j/2}\rho(2^j) \qquad r\ge 0.
$$
Hence, conditions on $\{\rho(n)\}$ will allow us to control
$\|\E[S_{{n}}|\F_0]\|_2$.

\medskip

We obtain
\begin{prop}\label{prop3}
Let $\{X_n\}$ be a stationary process such that $\rho(n)=O(\frac{1}
{(\log n)^{2}(\log \log n)^\tau}$ for some $\tau>1/2$, then \eqref{ZWC} holds.
\end{prop}
\noindent {\bf Remark:}\\
The proposition is based on Theorem \ref{lil}, while Theorem \ref{LIL}
does not really apply, since it seems that the only way to use $\{\rho(n)\}$
to check condition \eqref{WC} is via Proposition \ref{propC} and the estimate
 $\|\E[X_n|\F_0]\|\le \rho(n)$, which is not efficient.
{\bf Proof:}\\
Let $n\ge 1$ and $2^r\le n<2^{r+1}$. We have, for every $1\le l\le r$,
$$
\|\E[S_{2^{l}}|\F_0]\|_2\le C\sum_{j=0}^{l-1} \frac{2^{j/2}}{j^{2}
(\log j)^\delta} \le C' \frac{2^{l/2}}{l^{2}
(\log l)^\delta}.
$$
Hence,
$$
\|\E[S_{{n}}|\F_0]\|_2\le \sum_{l=0}^r\|\E[S_{{2^l}}|\F_0]\|_2
\le K \frac{\sqrt n }{(\log n)^{2} (\log \log n)^\delta}.
$$
\hfill $\square$

\begin{appendix}

\section{Proof of Proposition \ref{A1}}


\noindent Recall (see e.g. \cite{Zygmund} III, formula $(1.9)$), that, for
$0<\beta <1$,
\begin{equation}\label{zyg1}
(1-z)^{\beta-1}=\sum_{n\ge 0}a_n(\beta)z^n \qquad \forall |z|<1,
\end{equation}
where $\{a_n\}$ is decreasing to $0$ and the right hand side is well defined
for $|z|=1$, $z\neq 1$, and defines a continuous function on $\overline D
-\{1\}$. Moreover  (see e.g. \cite{Zygmund}, III.1, formula (1.18)),
\begin{equation}\label{Zyg2}
a_n(\beta)\underset{+\infty}\sim \frac{n^{-\beta}}{\Gamma(1-\beta)}(1+O(
\frac{1}{n})),
\end{equation}
where $\Gamma$ is the Euler's function.
We have

{\bf Proof:}\\
The convergence of the power series follows by Abel summation and the
non increasingness assumption on $b$.\\
By \eqref{zyg1} and \eqref{Zyg2}, we have, on $\overline D-\{1\}$
\begin{equation}\label{eq}
\sum_{n\ge 1}\frac{z^n}{n^\beta}= \Gamma(1-\beta)(1-z)^{\beta-1}
+o(|1-z|^{\beta-1}) \qquad z\rightarrow 1 .
\end{equation}

Let $\varepsilon >0$. Let $0<\omega <\Omega <+\infty$ and $z\in \overline
D-\{1\}$. Write
$$
\sum_{n\ge 1}\frac{z^n}{n^\beta} = \sum_{1\le n < \frac{\omega}{|1-z|}}
\frac{z^n}{n^\beta} +\sum_{\frac{\omega}{|1-z|}\le n\le \frac{\Omega}{|1-z|}}
\frac{z^n}{n^\beta} +\sum_ {n > \frac{\Omega}{|1-z|}}
\frac{z^n}{n^\beta}=S_1+S_2+S_3,
$$
and the corresponding decomposition
$$
\sum_{n\ge 1}\frac{b(n)z^n}{n^\beta}=T_1+T_2+T_3.
$$
We have
$$
|S_1|\le  \sum_{1\le n < \frac{\omega}{|1-z|}}\frac{1}{n^\beta} <
\int_0^{\frac{\omega}{|1-z|}}\frac{dt}{t^\beta}<\frac{\omega^{1-\beta}}
{(1-\beta)|1-z|^{1-\beta}}
$$
and, by Abel summation and the fact that $|\sum_{k=0}^nz^k|=\big|
\frac{1-z^{n+1}}{1-z}\big|\le \frac{2}{|1-z|}$,
$$
|S_3|\le \frac{2}{|1-z|}\big(\frac{|1-z|}{\Omega}\big)^{\beta}.
$$
So, one can find $\omega$ small enough and $\Omega$ large enough such that
$|S_1|\le \varepsilon |1-z|^{\beta-1}$, $|S_3| \le \varepsilon
|1-z|^{\beta-1}$. \\
Hence, by \eqref{eq}, for $z\in \overline D-\{1\}$ close to $1$, we have
$$
|S_2-\Gamma(1-\beta)(1-z)^{\beta-1}|\le 3\varepsilon |1-z|^{\beta-1}.
$$
We now estimate $T_1$ and  $T_3$ in the same way than $S_1$ and $S_3$ and
we will estimate $T_2$ thanks to $S_2$.\\
Let $\delta>0$ be fixed such that $\beta+\delta<1$. Since
$u^\delta b(u)$ is increasing to infinity,  there exists
$C\ge 1$ and $v>0$, such that for every $w\ge v$ and $1\le n < w$,
$n^\delta b(n)\le C w^\delta b(w)$. Hence,
provided that $\omega$ and $z$ are chosen such that
$ \frac{\omega}{|1-z|}\ge v$, we have
\begin{gather*}
|T_1|=|\sum_{1\le n < \frac{\omega}{|1-z|}}\frac{b(n)z^n}{n^\beta}|
\le \sum_{1\le n < \frac{\omega}{|1-z|}}\frac{b(n)n^\delta }
{n^{\beta+\delta}}\\
\le b(\frac{\omega}{|1-z|})\omega^\delta |1-z|^{-\delta}
\int_0^{\frac{\omega}{|1-z|}}
\frac{dt}
{t^{\beta+\delta}}\le \frac{\omega^{1-\beta}}{ 1-\beta-\delta}
b(\frac{\omega}{|1-z|})|1-z|^{\beta-1}.
\end{gather*}
Fix $\omega$ small enough such that $ \frac{\omega^{1-\beta}}{ 1-\beta-\delta}
\le \varepsilon$. Then, using that $b$ is slowly varying, for every $z$ close
enough to $1$, we have
$$
|T_1|\le 2\varepsilon b(\frac{1}{|1-z|})|1-z|^{\beta-1}.
$$

\noindent

If $\Omega/|1-z|$ is large enough, $u^{-\beta}b(u)$ is decreasing for $u\ge
\Omega/|1-z|$. Then we have, by Abel summation,
$$
|T_3|\le \frac{2}{|1-z|}b(\frac{1}{|1-z|})\big(\frac{|1-z|}{\Omega}\big)^\beta.
$$
So, for $\Omega$ large enough,
$$
|T_3|\le\varepsilon b(\frac{1}{|1-z|})|1-z|^{\beta -1}.
$$

\noindent

On the other hand, we have
$$
T_2=b(\frac{1}{|1-z|})\sum_{\frac{\omega}{|1-z|}\le n\le \frac{\Omega}{|1-z|}}
 \frac{z^n}{n^\beta}+ \sum_{\frac{\omega}{|1-z|}\le n\le \frac{\Omega}{|1-z|}}
 (b(n)-b(\frac{1}{|1-z|}))\frac{z^n}{n^\beta}=b(\frac{1}{|1-z|})S_2+T_2'.
$$
We already estimated $S_2$, and we have
$$
|T_2'|\le \max_{\frac{\omega}{|1-z|}\le n\le \frac{\Omega}{|1-z|}}|b(n)-
b(\frac{1}{|1-z|})|\sum_{\frac{\omega}{|1-z|}\le n\le \frac{\Omega}{|1-z|}}
\frac{1}{n^\beta}\le  \max_{\frac{\omega}{|1-z|}\le n\le
 \frac{\Omega}{|1-z|}}|b(n)-b(\frac{1}{|1-z|})| \frac{|1-z|^{\beta-1}}
{(1-\beta)}\Omega^{1-\beta}.
$$
Since $\omega$ and $\Omega$ are fixed and $b$ is slowly varying,
for $z$ close enough to $1$ we have
$|T_2'|\le \varepsilon b(\frac{1}{|1-z|})|1-z|^{\beta-1}$. \hfill $\square$
\\

\noindent {\bf Proof of Proposition \ref{A1}.}\\




We first estimate $B'$ and, by integration,  we will estimate $B$.

We have
$$
B'(z)=\sum_{n\ge 1}(\sum_{k\ge n} (\frac{c}{\sqrt{k^3b(k)}}) z^{n-1}
= \frac{c}{z-1}\sum_{k\ge 1}\frac{z^k-1}{\sqrt{k^3b(k)}},
$$
where the permutation of the sums is clearly justified.\\

We will write $B'(z):= c\frac{C(z)}{z-1}$. Then $C$ is analytic in
$D$ and continuous on $\overline D$ and $C(1)=0$. \\

For every $z\in D$, $C'(z)=\sum_{k\ge 1} \frac{z^{k-1}}{\sqrt{kb(k)}}$.
As previously, the series giving $C'$ defines a continous function
on $\overline D-\{1\}$.

Let $z\in \overline D-\{1\}$. Define $h_z(u)= C(1+u(z-1))$, $u\in [0,1].$
Then $h_z$ is differentiable on $]0,1]$, and for every $u\in ]0,1]$, we have
$$
h_z'(u)=(z-1)C'(1+u(z-1))=(z-1)\sum_{k\ge 1}\frac{(1+u(z-1))^{k-1}}
{\sqrt{kb(k)}}.
$$


Let $\varepsilon >0$. By Proposition \ref{prop1}, there exists $\delta>0$,
such that, for every $z\in \overline D-\{1\}$, with $|1-z|<\delta$ and every
$u\in ]0,1]$, we have, using $\Gamma(1/2)=\sqrt \pi$
$$
|\sum_{k\ge 1}\frac{(1+u(z-1))^{k-1}}{\sqrt{kb(k)}}-\frac{\sqrt \pi}
{\sqrt{b(\frac{1}{u|1-z|})}\sqrt{u(1-z)}}|\le \frac{\varepsilon}
{\sqrt{b(\frac{1}{u|1-z|})}\sqrt{u|1-z|}}.
$$
Hence
$$
\big|h_z'(u)-\frac{-\sqrt \pi \sqrt{1-z}}{\sqrt{ub(\frac{1}{u|1-z|})}}\big|
\le\frac{\varepsilon \sqrt{|1-z|}}{\sqrt{ub(\frac{1}{u|1-z|})}}.
$$
So, $h_z'$ is integrable and we have $C(z)-C(1)=\int_0^1h_z'(u)du$.
Hence, for
every $z\in \overline D-\{1\}$, with $|1-z|<\delta$,
\begin{gather*}
\bigg|C(z)-\int_0^1\frac{- \sqrt \pi \sqrt{1-z}}{\sqrt{ub(
\frac{1}{u|1-z|})}}du\bigg|
\le \varepsilon \sqrt{|1-z|}\int_0^1 \frac{du}{\sqrt{ub(\frac{1}{u|1-z|})}}\\
\le \varepsilon \int_0^{|1-z|} \frac{du}{\sqrt{ub(\frac{1}{u})}}
\le \frac{K_1\varepsilon \sqrt{|1-z|}}{\sqrt{b(\frac{1}{|1-z|})}},
\end{gather*}
where we used that $b$ is slowly varying.
Finally, we obtain
$$
\bigg|C(z)-\frac{- \sqrt \pi \sqrt{1-z}}
{\sqrt{|1-z|}}\int_0^{|1-z|}\frac{du}{\sqrt{ub(1/u)}}\bigg|
\le \frac{K_1\varepsilon \sqrt{|1-z|}}{\sqrt{b(\frac{1}{|1-z|})}}.
$$
Since, $b$ is slowly varying, we have (see e.g. \cite[Theorem 1.b]{Feller}
using the change of variable $t:=1/u$)
\begin{equation}\label{equiv}
\int_0^{x}\frac{du}{\sqrt {ub(1/u)}}
\underset{0}\sim \frac{2\sqrt x}{b(1/x)}.
\end{equation}
\medskip


For every $z\in \overline D-\{1\}$, define $g_z(v):= B(1+v(z-1))$,
$0\le v\le 1$.
Then $g_z$ is differentiable on $]0,1]$, and for every $v\in ]0,1]$, we
have
$$
g_z'(v)=(z-1)B'(1+v(z-1))=\frac{c}{v}C(1+v(z-1)).
$$
Hence, for every $z\in \overline D-\{1\}$, with $|1-z|<\delta$ and every $0<
v\le 1$, using \eqref{equiv}, we obtain
$$
\bigg|g_z'(v)-\frac{-2\sqrt \pi \sqrt{1-z}}{ |1-z|\sqrt v\sqrt{b(\frac{1}
{v|1-z|})}}\bigg| \le \frac{K\varepsilon \sqrt {|1-z|}}{\sqrt {vb(\frac{1}
{v|1-z|})}}.
$$
Thus $g'_z$ is integrable and $B(1)-B(z)=\int_0^1 g_z'(v)dv$. \\
Finally, using \eqref{equiv} and previous computations,
for every $z\in \overline D-\{1\}$, with $|1-z|<\eta$,
$|1-B(z)- \frac{-4\sqrt \pi \sqrt{1-z}}{b(\frac{1}{|1-z|})}|\le \frac{K'
\varepsilon \sqrt{|1-z|}}{b(\frac{1}{|1-z})}$, which proves the desired
result.
\bigskip

Let  prove $(iii)$ and $(ii)$.
\medskip


Recall that
$$
A(z)=\sum_{n\ge 0} \alpha_n z^n \qquad \forall z\in D.
$$
Moreover for every $0<r<1$ and $n\ge 0$,
\begin{equation}\label{coef}
\alpha_n=\frac{1}{2\pi r^n}\int_0^{2\pi} A(r{\rm e}^{it}){\rm e}^{-int}dt.
\end{equation}
Now observe that, for every $1/2\le r\le 1$ and
$0\le |t|\le \pi/2$, $|1-r{\rm e}^{it}|\ge r|\sin t|\ge \frac{2|t|r}{\pi}\ge
|t|/2$. Hence, by $(i)$, there exists $K>0$, such that, for every
$z=r{\rm e}^{it}$ close enough to $1$,
$$
|A(z)|\le \frac{K b(1/t)}{\sqrt t},
$$
which defines an integrable function with respect to $t$.
Since $A$ is continous on $\overline D-\{1\}$, changing $K$ if necessary, we
have
$$
\sup_{1/2\le r\le 1} |A(r{\rm e}^{it}|\le \frac{K b(1/t)}{\sqrt t}.
$$
Thus, by Lebesgue dominated convergence theorem and \eqref{coef}
\begin{equation}\label{coef2}
\alpha_n=\frac{1}{2\pi}\int_0^{2\pi} A({\rm e}^{it}){\rm e}^{-int}dt \qquad
\forall n\ge 0.
\end{equation}

Similarly,
$$
0=\lim_{r\rightarrow 1^-}\int_0^{2\pi}A(r{\rm e}^{it}){\rm e}^{int}dt
=\int_0^{2\pi}A({\rm e}^{it}){\rm e}^{int}dt.
$$
Denote $a_1(t):=A({\rm e}^{it})$, $0< t<2\pi$. Then $a_1$ is integrable and
admits $\{\alpha_n\}_{n\ge 0}$ for Fourier coefficients. Since, by Corollary 3
of \cite{ZW}, $a_1$ is continuously differentiable, hence its Fourier series
converging to $a_1$ on $]0,2\pi [$. \\
Hence $(iii)$ is proven.\\


To prove $(ii)$, we use  Proposition 3 of \cite{ZW}: we have
$|\alpha_n-\alpha_{n+1}|=O(\frac{\sqrt {b(n)}}{\sqrt {n^3}})$. Since
$\{\alpha_n\}$ are the Fourier coefficients of an integrable function,
$\alpha_n\rightarrow 0$, and
\begin{equation}\label{a_n}
|\alpha_n|\le \sum_{k\ge n}|\alpha_k-\alpha_{k+1}|=O(\frac{\sqrt {b(n)}}
{\sqrt n}).
\end{equation}
Let $z\in \overline D-\{1\}$ and $m\ge 0$. \\
If $\frac{1}{|1-z|}\ge m$, then
\begin{gather*}
|\sum_{n=0}^m \alpha_nz^n|\le \sum_{n=0}^m |\alpha_n|
\le K \sqrt {mb(m)} \le K \frac{\sqrt{b(\frac{1}{|1-z|})}}{\sqrt{|1-z|}},
\end{gather*}
which proves the result in that case.\\
Assume now that $\frac{1}{|1-z|}\le m$.
Let $S_n=\sum_{k=0}^n z^k$,
$n\ge 0$. We have
\begin{align*}
\sum_{n=0}^m \alpha_nz^n  &  = A(z)-\sum_{n\ge m+1} \alpha_n(S_n-S_{n-1})\\
                          &  =A(z)-\sum_{n\ge m+1}S_n (\alpha_n-\alpha_{n+1}) +
\alpha_{m+1}S_m.
\end{align*}
We already saw that $|S_n|\le \frac{2}{|1-z|}$. Hence by \eqref{a_n} and
$(i)$, we have
\begin{equation}
|\sum_{n=0}^m \alpha_nz^n |\le |A(z)|+\frac{2}{|1-z|}
\frac{K\sqrt {b(m)}}{\sqrt m})\le \tilde K \frac{\sqrt{b(\frac{1}{|1-z|})}}
{\sqrt{|1-z|}},
\end{equation}
which gives the result in that case. \hfill $\square$
\\


\section{Proof of Proposition \ref{prop}}

Recall that for every $0\le t<1$, $\Gamma_t:=\sum_{n\ge 0}t^nX_n$ and that
$\varphi_t:= \E[\Gamma_t\circ \theta |\F_1]-\E[\Gamma_t\circ \theta |\F_0]$.\\
By \eqref{Abel} and our assumption, there exists $K>0$ such that
\begin{gather}
\|\E[\Gamma_t|\F_0]\|_2\le K(1-t)\sum_{n\ge 0}t^n \frac{\sqrt n}{(\log (n+1))^2
(\log \log (n+2))^\tau} \nonumber\\
\label{Gamma}=O \big(\frac{1}{\sqrt {1-t} (\log(1-t))^2
(\log|\log (1-t)|)^\tau}\big) \quad (t\rightarrow 1),
\end{gather}
by a Tauberian theorem (see Theorem 5 in  \S XIII.5 of \cite{Feller}).
A similar estimate can be obtained for $\|\E[\Gamma_t\circ \theta|\F_0]\|_2$.\\

Let $0<t <s<1$. We have that for every $U\in L^2$,
$$\E\big[((\E[U|\F_0])\circ \theta-\E[U\circ \theta|\F_0])^2\big]=
\E[((\E[U|\F_0])^2]-\E[((\E[U\circ \theta |\F_0])^2].$$
Taking $U:=\Gamma_t-\Gamma_s$, we obtain
\begin{align*}
\E[(\varphi_s-\varphi_t)^2]  &  =\E\big[((\E[U|\F_0])^2]-
\E[((\E[U\circ \theta |\F_0])^2\big]\\
  &  =\E \big[(\E[U|\F_0]-\E[U\circ \theta |\F_0])
(\E[U|\F_0]+\E[U\circ \theta |\F_0])\big]\\
  &  \le \|\E[U|\F_0]-\E[U\circ \theta |\F_0]\|_2
 \cdot  \|\E[U|\F_0]+\E[U\circ \theta |\F_0]\|_2.
\end{align*}
Since $t\E[G_t\circ \theta |\F_0]=\E[G_t|\F_0]-X_0$, we have for $t>1/2$,
\begin{gather*}
\E[(\varphi_s-\varphi_t)^2] \le
\big(2(1-t)\|\E[\Gamma_t|\F_0]\|_2+2(1-s)
\|\E[\Gamma_s|\F_0]\|_2+(s-t)\|X_0\|_2\big) \\
\big(\|\E[\Gamma_t|\F_0]\|_2
+\|\E[\Gamma_t\circ \theta|\F_0]\|_2+\|\E[\Gamma_s|\F_0]\|_2
+\|\E[\Gamma_s\circ \theta|\F_0]\|_2\big)
\end{gather*}
Hence, using \eqref{Gamma}, for every $t>1/2$ and every $t<s<(1+t)/2$, we have
\begin{gather*}
\E[(\varphi_s-\varphi_t)^2] \le K'(\frac{\sqrt{1-t}}{ (\log(1-t))^2
(\log|\log (1-t)|)^\tau}+\frac{\sqrt{1-s}}{ (\log(1-s))^2
(\log|\log (1-s)|)^\tau}+s-t)\\
( \frac{1}{\sqrt {1-t} (\log(1-t))^2
(\log|\log (1-t)|)^\tau} +\frac{1}{\sqrt {1-s} (\log(1-s))^2
(\log|\log (1-s)|)^\tau})\\
\le\frac{K''}{(\log (1-t))^4(\log |\log (1-t)|)^{2\tau}}.
\end{gather*}
For every $n\ge 0$, write $t_n:=(2^n-1+t)/2^n$. Then
$$
\sup_{t_n\le s<t_{n+1}}\|\varphi_s-\varphi_{t_n}\|_2\le
\frac{\sqrt{K''}}{(|\log (1-t)|+n\log 2)^2
(\log (|\log (1-t)|+n\log 2))^{\tau}}.
$$
Summing over $n$ we obtain, for every $0\le t<1$,
\begin{gather*}
\sup_{t\le s<1}\|\varphi_s-\varphi_{t}\|_2\le \sqrt{K''}\int_0^\infty
\frac{dx}{(|\log (1-t)|+x\log 2)^2(\log (|\log(1-t)|+x\log 2))^\tau}\\
=\frac{\sqrt{K''}}{\log 2}\int_{|\log (1-t)|}^\infty \frac{du}
{u^2(\log u)^\tau}\le \frac{\tilde K}{|\log (1-t)|(\log |\log (1-t)|)^\tau}.
\end{gather*}
By Cauchy's criterion, there exists $M\in L^2(\Omega,\P)$, such that
$\{\varphi_t\}$ converges in $L^2$ to $M$. Moreover
$$
\|M-\varphi_t\|_2\le \frac{\tilde K}{|\log (1-t)|(\log |\log (1-t)|)^\tau}
\qquad \forall t\in [0,1).
$$
Define $M_n:=\sum_{k=1}^n M\circ \theta^k$. Then $\{M_n\}$ is a martingale
with
stationary increments and
$$
\|M_n-M_n(t)\|_2\le \frac{\sqrt n\tilde K}{|\log (1-t)|
(\log |\log (1-t)|)^\tau}\qquad \forall t\in [0,1).
$$
Now
\begin{gather*}
S_n-M_n= \sum_{k=1}^n (\E[\Gamma_t|\F_0]-t\E[\Gamma_t\circ \theta|\F_0]
)\circ \theta^k-M_n(t)+M_n(t)-M_n\\
=\E[\Gamma_t|\F_0]\circ \theta -\E[\Gamma_t|\F_0]\circ \theta^{n+1}
+(1-t)\sum_{k=1}^n \E[\Gamma_t\circ \theta|\F_0]\circ \theta^k
+M_n(t)-M_n.\\
\end{gather*}
Hence, using \eqref{X_0}, we obtain for every $t\in [0,1)$, $n\ge 1$,
\begin{gather*}
 \|S_n-M_n\|_2\le 2\|\E[\Gamma_t|\F_0]\|_2+n(1-t)\|\E[\Gamma_t\circ\theta
 |\F_0]\|_2+\|M_n(t)-M_n\|_2
\end{gather*}
Taking $t=1-1/n$ in the above and using \eqref{Gamma} (see also the remark
after it), we deduce
$$
\|S_n-M_n\|_2 \le  \frac{K\sqrt n}{\log n (\log \log n)^\tau},
$$
which finishes the proof. \hfill $\square$

\section{Proof of Proposition \ref{normal}}

The proof follows the lines of the proof of Proposition \ref{prop}, except
that we will use the spectral calculus to obtain better estimates. The
beginning of the proof is closely related to the work of
Kipnis-Varadha \cite{KV} for symmetric Markov chains.\\
For every $0\le t<1$ define, $G_tf:=\sum_{n\ge 0} t^nP^nf$. Then $G_tf(X_0)$
corresponds to $\E[\Gamma_t|\F_0]$, in the proof of Proposition \ref{prop}.
Define $\varphi_t:=G_tf(X_1)-PG_tf(X_0)$.

We would like to prove that $\{\varphi_t\}$ converges in $L^2(\P_m)$
when $t\rightarrow 1$.\\
Let $0\le s <t <1$. We have
\begin{gather*}
\E_m[(\varphi_t-\varphi_s)^2]= \E_m[(G_t-G_sf)^2(X_0)]-
\E_m[(PG_t-PG_sf)^2(X_0)].
\end{gather*}
Let $\mu_f$ be the spectral measure of $f$ with respect to the normal operator
$P$. It is not difficult to see that
\begin{gather*}
\E_m[(G_t-G_sf)^2(X_0)]-\E_m[(PG_t-PG_sf)^2(X_0)]\\=
\int_D\big|\frac{1}{1-tz}-\frac{1}{1-sz}\big|^2\mu_f(dz)-
\int_D\big|\frac{z}{1-tz}-\frac{z}{1-sz}\big|^2\mu_f(dz)\\
=(1-|z|^2)\int_D\big|\frac{1}{1-tz}-\frac{1}{1-sz}\big|^2\mu_f(dz)\\
\le 2 \int_D\big|\frac{|1-z|^{1/2}}{1-tz}-\frac{|1-z|^{1/2}}
{1-sz}\big|^2\mu_f(dz).
\end{gather*}
Hence it is enough to show that when $t\rightarrow 1$, $\displaystyle
\{\frac{|1-z|^{1/2}}{1-tz}\}_t$ converges in $L^2(\mu_f)$.\\
Let $z:=r{\rm e}^{i\theta}$, with $\cos \theta\ge 0$.\\
If $r\le \cos \theta$, one can see geometrically that
$$
|1-tz|\ge |1-z| \qquad \forall 0\le t \le 1.
$$
Assume now that $r>\cos \theta$. One can see geometrically that
$|1-tz|^2\ge \sin^2\theta =1-\cos^2\theta \ge 1-\cos \theta$. On the other
hand simple computations show that, since  $r>\cos \theta$, $|1-z|^2\le
|1-{\rm e}^{i\theta} |^2=2(1-\cos\theta)$. \\
Finally, for every $z=r{\rm e}^{i\theta}$, with $\cos \theta \ge 0$,
$\sqrt 2|1-tz|\ge |1-z|$ for every $t\in [0,1]$. Hence there exists
$K>0$ such that
$$
\frac{|1-z|^{1/2}}{|1-tz|}\le \frac{K}{|1-z|^{1/2}}
 \qquad \forall z\in D,~\forall t\in [0,1].
$$
Since, by Proposition \ref{sqrt}, $\int_D\frac{\mu_f(dz)}{|1-z|}<+\infty$,
Lebesgue dominated theorem yields to the convergence in $L^2(\mu_f)$ of
$\displaystyle \{\frac{|1-z|^{1/2}}{1-tz}\}_t$ to $\frac{|1-z|^{1/2}}{1-z}$. \\
Hence $\{\varphi_t\}$ is a Cauchy sequence in $L^2(\Omega,\P_m)$, so it
converges
 to an element $M\in L^2[\Omega,\P_m)$. Moreover, for every $0\le t<1$, we have
$$
\E_m[(\varphi_t-M)^2]\le 2(1-t)^2\int_D \frac{\mu_f(dz)}{|1-z||1-tz|^2}.
$$
Define, for every $n\ge 1$, $0\le t <1$
$$
M_n(t):= \sum_{k=0}^{n-1} \varphi_t \circ \theta^k \mbox{ and }
M_n:=\sum_{k=0}^{n-1} M\circ \theta^k.
$$
Then $\{M_n(t)-M_n\}$ is a martingale with stationary increments and
$$
\E_m[(M_n(t)-M_n)^2]=n\E_m[(\varphi_t-M)^2]\le 2n (1-t)^2\int_D
\frac{\mu_f(dz)}{|1-z||1-tz|^2}.
$$
Write $W_n:=S_n(f)-M_n$. It remains to estimate $\E_m[W_n^2]$. We have
$$
W_n=S_n(f)-M_n(t)+M_n(t)-M_n=M_n(t)-M_n +G_tf(X_0)-G_tf(X_n)+(1-t)
\sum_{k=0}^{n-1}PG_tf(X_k).
$$
Hence, since, for evey $l\in\{0,1\}$ and $0\le k\le n-1$,
$\E_m[(P^lG_tf(X_k))^2]=\int_D \frac{|z|^{2l}}{|1-tz|^2}\mu_f(dz)$, we have
\begin{gather*}
\E_m[W_n^2]\le 4 (\E [( M_n(t)-M_n)^2]+\E_m[( G_tf(X_0))^2]+
\E_m[(G_tf(X_n))^2]+ (1-t)^2\E_m[(\sum_{k=0}^{n-1}PG_tf(X_k))^2]\\
\le 4\big(2n (1-t)^2\int_D \frac{\mu_f(dz)}{|1-z||1-tz|^2} +2 \int_{D_1}
\frac{\mu_f(dz)}{|1-tz|^2}\mu_f(dz)+n^2(1-t)^2\int_{D_1}
\frac{\mu_f(dz)}{|1-tz|^2}\big).
\end{gather*}
Take $t=u_n:=1-1/n$ to obtain
$$
\E_m[W_n^2]\le 8\big(\frac{1}{n}\int_{D_1} \frac{\mu_f(dz)}{|1-z||1-u_nz|^2}+
\int_{D_1} \frac{\mu_f(dz)}{|1-u_n z|^2}\big).
$$
We already saw that there exists $C>0$ such that $C |1-u_nz|\ge |1-z|$ for
every $z\in D$ and every $n\ge 1$. Moreover we also have $|1-u_nz|\ge
|1-u_n|=1/n$. \\
Now recall that $D_n=\{z=r{\rm e}^{2i\pi \theta}~:~1-\frac{1}{n}\le r
\le 1,~-\frac{1}{n}\le \theta \le \frac{1}{n}\}$.\\
We have
\begin{gather*}
\E_m[W_n^2]\le C_1 \bigg(\frac{1}{n}\big( \int_{D_1-D_n}\frac{\mu_f(dz)}
{|1-z|^3}+n^2\int_{D_n}\frac{\mu_f(dz)}{|1-z|}\big) +
\big( \int_{D_1-D_n} \frac{\mu_f(dz)}{|1-z|^2}+n^2\mu_f(D_n)\bigg)\\
\le C_3\big( \int_{D_1-D_n} \frac{\mu_f(dz)}{|1-z|^2}+
n\int_{D_n}\frac{\mu_f(dz)}{|1-z|}\big).
\end{gather*}
Now
\begin{gather*}
 \int_{D_1-D_n} \frac{\mu_f(dz)}{|1-z|^2} =\sum_{k=1}^{n_1}
\int_{D_k-D_{k+1}} \frac{\mu_f(dz)}{|1-z|^2}\\
\le C_4 \sum_{k=1}^{n-1} k^2(\mu_f(D_k)-\mu_f(D_{k+1}))
= C_4 \sum_{k=1}^{n-1} \mu_f(D_k)(k^2-(k-1)^2) -(n-1)^2\mu_f(D_n)\\
\le 2C_4 \sum_{k=1}^{n-1} k \mu_f(D_k)
\le C_5 \sum_{k=1}^{n-1} \frac{\|U_k(f)\|_2^2}{k},
\end{gather*}
where the last inequality comes from \eqref{Sn}.\\
On the other hand, we have
\begin{gather*}
n\int_{D_n}\frac{\mu_f(dz)}{|1-z|}=n\sum_{k\ge n}\int_{D_k-D_{k+1}}
\frac{\mu_f(dz)}{|1-z|}\\
\le nC_6 \sum_{k\ge n}k(\mu_f(D_k)-\mu_f(D_{k+1}))
=nC_6 (n\mu_f(D_n)+\sum_{k\ge n+1 }\mu_f(D_k))\\
\le nC_6\big( \frac{\|U_n(f)\|_2^2}{n}+\sum_{k\ge n+1} \frac{\|U_k(f)\|_2^2 }
{k^2}\big),
\end{gather*}
which finishes the proof of the proposition. \hfill $\square$

\noindent {\bf Acknowledgement.} The author is very grateful to Guy Cohen for
 a careful reading of the manuscript, and to Michael Lin and Guy Cohen for
valuable discussion.

\end{appendix}

\end{document}